\documentclass[12pt]{article}
\usepackage{amsmath,amssymb,graphicx,curves}
\newcommand{\qed}{\hfill\vrule height6pt  width6pt depth0pt \medskip}
\newcommand{\ppt}{\textrm{PPT}}
\newcommand{\sppt}{\textrm{SPPT}}
\newcommand{\qpp}{{\cal Q}}
\newcommand{\qppc}{\overline{\cal Q}}

\newcounter{mycount}
\newenvironment{mylist}{\begin{list}{(\roman{mycount})}%
{\usecounter{mycount}\itemsep 0pt}}{\end{list}}
\title{The dynamics of Pythagorean triples}
\author{Dan Romik
\footnote{
Department of Mathematics, Weizmann Institute of Science,
Rehovot 76100, Israel. email: \texttt{romik@wisdom.weizmann.ac.il}}}

\begin{document}
\maketitle

\begin{abstract}
    We construct a piecewise onto 3-to-1 dynamical system on the positive
    quadrant of the unit circle, such that for rational points (which
    correspond to normalized Primitive Pythagorean Triples), the associated
    ternary expansion is finite, and is equal to the address of the PPT on
    Barning's \cite{barning} ternary tree of PPTs, while irrational points
    have infinite expansions. The dynamical system is conjugate to a
    modified Euclidean algorithm. The invariant measure is identified, and
    the system is shown to be conservative and ergodic. We also show, based
    on a result of Aaronson and Denker \cite{aarodenker}, that the
    dynamical system can be obtained as a factor map of a cross-section of
    the geodesic flow on a quotient space of the hyperbolic plane by the
    group $\Gamma(2)$, a free subgroup of the modular group with two
    generators.
\end{abstract}

\section{Introduction}

The starting point of this paper is a theorem, attributed to Barning
\cite{barning}, on the structure of the set of Primitive Pythagorean
Triples, or PPTs. Recall that a PPT is a triple $(a,b,c)$ of integers,
with $a,b,c>0, \textrm{ gcd}(a,b)=1$ and
\begin{equation}\label{eq:pythagoras}
a^2+b^2=c^2.
\end{equation}
Clearly, if $(a,b,c)$ is a PPT, then one of $a,b$ must be odd, and the
other even. Barning \cite{barning}, and later independently several
others \cite{alperin, gollnick, hall, jaeger, kanga, kristensen,
preau} (see also \cite{mccullough}), showed:

{\it \paragraph{Theorem 1.} Define the matrices
\begin{equation}\label{eq:matrices}
M_1 = \left( \begin{array}{rrr} -1 & 2 & 2 \\
                          -2 & 1 & 2 \\
                          -2 & 2 & 3
       \end{array} \right),\ \ 
M_2 = \left( \begin{array}{rrr} 1 & 2 & 2 \\
                          2 & 1 & 2 \\
                          2 & 2 & 3
       \end{array} \right),\ \ 
M_3 = \left( \begin{array}{rrr} 1 & -2 & 2 \\
                          2 & -1 & 2 \\
                          2 & -2 & 3
       \end{array} \right).
\end{equation}
Any PPT $(a,b,c)$ with $a$ odd and $b$ even has a unique
representation as the matrix product
\begin{equation}\label{eq:expansion1}
\left( \begin{array}{l}a\\b\\c \end{array}\right) =
M_{d_1}M_{d_2}\ldots M_{d_n}
\left(\begin{array}{l}3\\4\\5\end{array}\right),
\end{equation}
for some $n\ge 0$,\ \  $(d_1,d_2,\ldots,d_n)\in\{1,2,3\}^n$.
Any PPT $(a,b,c)$ with $a$ even and $b$ odd has a unique
representation as
\begin{equation}\label{eq:expansion2}
\left( \begin{array}{l}a\\b\\c \end{array}\right) =
M_{d_1}M_{d_2}\ldots M_{d_n}
\left(\begin{array}{l}4\\3\\5\end{array}\right),
\end{equation}
for some $n\ge 0$,\ \  $(d_1,d_2,\ldots,d_n)\in\{1,2,3\}^n$. Any
triple of one of the forms \eqref{eq:expansion1},
\eqref{eq:expansion2} is a PPT.
}

\bigskip
In some of the papers where this was discussed, the theorem has been
described as placing the PPTs $(a,b,c)$ with $a$ odd and $b$ even on
the nodes of an infinite rooted ternary tree, with the root
representing the ``basic'' triple $(3,4,5)$, and where each triple
$(a,b,c)$ has three children, representing the multiplication of the
triple (considered as a column vector) by the three matrices $M_1,
M_2, M_3$. In this paper, we consider a slightly different outlook. We
think of the sequence $(d_1,d_2,\ldots,d_n)$ in \eqref{eq:expansion1},
\eqref{eq:expansion2} as an \emph{expansion} corresponding to the
triple $(a,b,c)$, over the ternary alphabet $\{1,2,3\}$. We call the
$d_i$'s the \emph{digits} of the expansion. To distinguish between
PPTs with $a$ odd, $b$ even and those with $a$ even, $b$ odd, we affix
to the expansion a final digit $d_{n+1}$, which can take the values
$oe$ ($a$ odd, $b$ even) or $eo$ ($a$ even, $b$ odd). So we have a 1-1
correspondence
$$ (a,b,c)\in \ppt\ \ \longleftrightarrow\ \ (d_1,d_2,\ldots,d_{n+1}) \in
\bigcup_{n=0}^\infty \{1,2,3\}^n \times \{ oe, eo \}.
$$
Several questions now come to mind:

\begin{itemize}
\item Is there a simple way to compute the expansion of a PPT? {\it
(Yes -- this is contained in the proof of Theorem 1.)}

\item As is easy to see and has been known since ancient times, the
mapping $(a,b,c)\to(a/c,b/c)$ gives a 1-1 correspondence between the
set of PPTs and the rational points $(x,y)$ on the positive quadrant
$\qpp$ of the unit circle. Are the digits in the expansion piecewise
continuous functions on the quadrant? {\it (Yes.)} Can one define a
ternary expansion for \emph{irrational points}? {\it (Yes.)}

\item Can interesting things be said about the expansion
$(d_1,d_2,\ldots,d_{n+1})$ of a \emph{random} PPT, chosen from some
natural model for random PPTs, say by choosing uniformly at random
from all PPTs $(a,b,c)$ with $c\le N$ and letting $N\to\infty$? {\it
(Yes.)}

\item Do these questions lead to interesting mathematics? {\it (Yes --
they lead to a dynamical system on $\qpp$ with interesting
properties.)}
\end{itemize}

It is the goal of this paper to answer these questions. The basic
observation is that it is in many ways preferable to deal with points
$(x,y)$ on the positive quadrant $\qpp$ of the unit circle, instead of
with PPTs. In the proof of Theorem 1, we shall see that there is a
simple transformation (a piecewise linear mapping) that takes a PPT
$(a,b,c)$ to the PPT $(a',b',c')$ that corresponds to its
\emph{parent} on the ternary tree, i.e., if $(a,b,c)$ has expansion
$(d_1,\ldots,d_{n+1})$ then $(a',b',c')$ has expansion
$(d_2,\ldots,d_{n+1})$. A standard trick in dynamical systems is to
rescale such transformations, discarding information that is
irrelevant for the continuing application of the transformation; this
is done, for example, when rescaling the Euclidean algorithm mapping
$(x,y)\to (y,x\textrm{ mod }y)$ to obtain the continued fraction
transformation $x\to\{1/x\}$. When we apply this idea to our case, we
obtain the following result.

{\it \paragraph{Theorem 2.} Let
$$ \qpp = \{ (x,y) : x>0,\ y>0,\ \ x^2 + y^2 = 1 \}. $$
Define the transformation $T:\qpp\to\qppc$ by
$$ T(x,y) =
   \left(\frac{|2-x-2y|}{3-2x-2y},\frac{|2-2x-y|}{3-2x-2y}\right).
$$
Define $d:\qpp\to\{1,2,3,oe,eo\}$ by
$$ d(x,y) = \left\{ \begin{array}{lllllll}
   1 & \quad & 4/3 &< & x/y, & & \\ %\frac{y}{x} > \frac{4}{3}, \\
   2 & & 3/4 &<& x/y& <& 4/3, \\ %\frac{3}{4} < \frac{y}{x} < \frac{4}{3}, \\
   3 & & & &x/y& <& 3/4, \\\ \\ %\frac{y}{x} < \frac{3}{4}, \\
   oe & & (x,y)& =& \left(\frac{3}{5},\frac{4}{5}\right), \\
   eo & & (x,y)& =& \left(\frac{4}{5},\frac{3}{5}\right).
\end{array}\right. $$
Then:
\begin{mylist}
\item If $(x,y)=(a/c,b/c)\in\qpp\cap\mathbb{Q}^2$ is a rational point of
$\qpp$, with $a/c,\ b/c$ in lowest terms (so $(a,b,c)$ is a PPT), then
for some $n\ge 0$, $T^{n+1}(x,y)$ (the $(n+1)$-th iterate of $T$) will
be equal to $(1,0)$ or $(0,1)$, and if we define
$$ d_k = d(T^{k-1}(x,y)), \qquad k=1,2,\ldots,n+1, $$
then $(d_1,d_2,\ldots,d_{n+1})$ is the ternary expansion (with the
last digit in $\{oe,eo\}$) corresponding to the PPT $(a,b,c)$ as in
\eqref{eq:expansion1}, \eqref{eq:expansion2}.

\item If $(x,y)\in\qpp$ is an irrational point, then $T^n(x,y)\in\qpp$
for all $n\ge 0$, and the sequence
$$ d_k = d(T^{k-1}(x,y)), \qquad k\ge 1, $$ 
defines an infinite expansion for $(x,y)$ over the alphabet
$\{1,2,3\}$, with the property that it does not terminate with an
infinite succession of $1$'s or with an infinite succession of $3$'s.

\item Any sequence $(d_k)_{k\ge 1}$ over the alphabet $\{1,2,3\}$
which does not terminate with an infinite succession of $1$'s or an
infinite succession of $3$'s determines a unique (irrational) point
$(x,y)\in\qpp$ such that $d_k = d(T^{k-1}(x,y)),\ \ k\ge 1$.
\end{mylist}
}

\paragraph{Examples.} Here are some examples of points in $\qpp$ and
their expansions. If $(d_k)_{1\le k\le n}$ is an expansion (finite or
infinite), we denote by $[d_1,d_2,\ldots]$ the point $(x,y)\in\qpp$
which has the given sequence as its expansion.
$$ \begin{array}{lll}
(3/5,4/5) &=& [oe]\\ (4/5,3/5)&=&[eo]
\end{array} \qquad
\begin{array}{lll}
(15/17,8/17) &=& [1,oe]\\ (21/29,20/29) &=& [2,oe]\\
(5/13,12/13) &=& [3,oe] 
\end{array}
$$ 
{\scriptsize $$
\begin{array}{lll}
(35/37,12/37)&=&[1,1,oe]\\
(77/85,36/85)&=&[1,2,oe]\\
(45/53,28/53)&=&[1,3,oe]
\end{array}\ \ 
\begin{array}{lll}
(65/97,72/97)&=&[2,1,oe]\\
(119/169,120/169)&=&[2,2,oe]\\
(55/73,48/73)&=&[2,3,oe]
\end{array}\ \ 
\begin{array}{lll}
(33/65,56/65)&=&[3,1,oe]\\
(39/89,80/89)&=&[3,2,oe]\\
(7/25,24/25)&=&[3,3,oe]
\end{array}
$$}
$$ \begin{array}{lll}
(\sqrt{2}/2,\sqrt{2}/2) &=&
[2,2,2,2,\ldots] \\
(1/2,\sqrt{3}/2) &=&
[3,1,3,1,\ldots] \\
(\sqrt{3}/2,1/2) &=&
[1,3,1,3,\ldots]
\end{array} \quad
\begin{array}{lll}(2/\sqrt{5},1/\sqrt{5})&=&[1,2,1,2,\ldots] \\
(1/\sqrt{5},2/\sqrt{5})&=&[3,2,3,2,\ldots] \\
(3/\sqrt{10},1/\sqrt{10})&=&[1,1,2,1,1,2,\ldots] \end{array}
$$ $$ \begin{array}{lll}
(\cos 1,\sin 1) &\ \ \, =& [3,1,1,3,1,1,1,1,3,1,1,1,1,1,1,3,1,1,1,1,1,
1,1,1,3,\ldots ]
\end{array}
\qquad\qquad\qquad\qquad\qquad\qquad\qquad\qquad
$$
\quad (see section 5)
$$
(\cos(1/\pi),\sin(1/\pi))\ \ =\  [
1, 1, 2, 1, 2, 2, 3, 3, 3, 3, 3, 3, 3, 1, 3, 3, 3, 3, 3, 3, 2,\ldots]
$$
\quad (this is meant as an example of a ``typical'' expansion
-- see section 4)

\vspace{-15.0pt}
\begin{eqnarray*}
[1,1,\ldots,1,oe]\ \ \textrm{($n$ times ``1'')} &=&
\left(\frac{4(n+1)^2-1}{4(n+1)^2+1},
\frac{4(n+1)}{4(n+1)^2+1}\right) \\
\ [2,2,\ldots,2,oe]\ \ \textrm{($n$ times ``2'')} &=&
\left(\frac{a_n}{c_n}, \frac{a_n+(-1)^n}{c_n}\right)
\end{eqnarray*}
where $(a_n)_{n\ge 0}=(3,21,119,697,\ldots)$ and $(c_n)_{n\ge 0}
=(5,29,169,985,\ldots)$ are sequences A046727 and A001653,
respectively, in The On-Line Encyclopedia of Integer Sequences
\cite{sloane}.

\bigskip
After constructing the dynamical system associated with the ternary
expansion, the next step is to study its properties. What does
the expansion of a typical point look like? To a trained eye, the
answer is contained in the following theorem.

{\it \paragraph{Theorem 3.} Let $ds$ denote arc length on the unit
circle. The dynamical system $(\qpp, T)$ possesses an infinite
invariant measure $\mu$, given by
$$ d\mu(x,y) = \frac{ds}{\sqrt{(1-x)(1-y)}}. $$
With the measure $\mu$, the system $(\qpp,T,\mu)$ is a conservative
and ergodic infinite measure-preserving system.
%The entropy of the system is
%$$ \kappa = \sqrt{2}\left( \frac{\pi^2}{12}+\frac{3(\log 3)^2}{2} -
%(\log 3)(\log 4)+\li(-1/3)+2 \li(2/3)\right),
%$$
%where $\li$ is the dilogarithm function
%$$ \li(x)=\sum_{n=1}^\infty \frac{x^n}{n^2}=
%\int_0^1-\frac{\log(1-x)}{x}dx.$$
}

\bigskip The invariant measure $\mu$ encodes all the information about the
statistical regularity of expansions of ``typical'' points.  In section 4
we shall state more explicitly some of the number-theoretic consequences of
Theorem 3.

Recall (\cite{hardywright}, Theorem 225) that the general parametric
solution of the equation \eqref{eq:pythagoras} with $a,b>0$ coprime,
$a$ odd and $b$ even is given by
\begin{equation}\label{eq:parametric}
a = m^2-n^2, \qquad b=2mn, \qquad c=m^2+n^2,
\end{equation}
where $m,n$ have opposite parity, $m>n>0$, and $\textrm{gcd}(m,n)=1$. A
roughly equivalent statement is that the map
$$ D:t\longrightarrow
\left(\frac{1-t^2}{1+t^2},\frac{2t}{1+t^2}\right)
$$
maps the extended real line injectively onto the unit circle, and maps
the rational numbers, together with the point at infinity, onto the
rational points of the circle. Note also that the interval $(0,1)$ is
mapped onto the positive quadrant $\qpp$. It is thus natural, in
trying to understand the behavior of the dynamical system $(\qpp,T)$,
to conjugate it by the mapping $D$, to obtain a new dynamical system
on $(0,1)$. This leads to the following result (the precise meaning of
the last statement in the theorem will be explained later):

{\it \paragraph{Theorem 4.} $(\qpp,T,\mu)$ is conjugate, by the
mapping $D$, to the measure preserving system $((0,1),\hat{T},\nu)$,
where
\begin{eqnarray*}
\hat{T}(t) &=& (D^{-1}\circ T\circ D)(t) = \left\{
  \begin{array}{ll} \frac{t}{1-2t} & 0<t<\frac{1}{3}, \\
    \frac{1}{t}-2 & \frac{1}{3}<t<\frac{1}{2}, \\
    2-\frac{1}{t} & \frac{1}{2}<t<1, \end{array}\right. \\
d\nu(t) &=& \frac{1}{\sqrt{2}}\cdot \frac{dt}{t(1-t)}.
\end{eqnarray*}
The mapping $\hat{T}$ is the scaling of a modified slow (subtractive)
Euclidean algorithm.
}

\begin{figure}[h!]
\begin{center}
\includegraphics{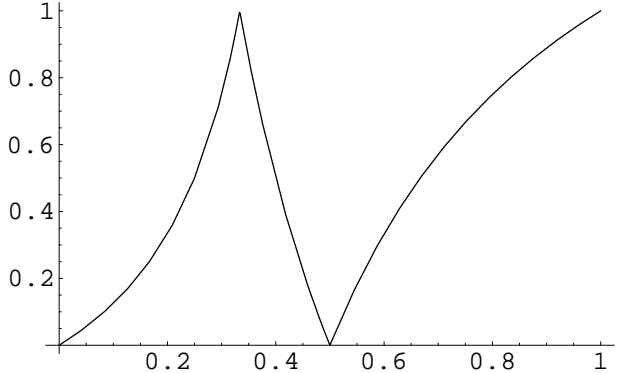}
\caption{The interval map $\hat{T}$}
\end{center}
\end{figure}

\bigskip
%The observation that $\hat{T}$, and therefore its conjugate map $T$,
%are related to subtractive Euclidean algorithms, can be used to
%accelerate the algorithm which by its nature belongs to the family of
%``slow'' algorithms. (In dynamical systems parlance, this is related
%to the fact that the map $\hat{T}$ has two \emph{indifferent fixed
%points}; that is also why the invariant measure is infinite.) This has
%certain computational advantages -- see section 4 for discussion.
%
There has been some interest in obtaining natural dynamical systems
appearing in number theory as factors of certain cross sections of the
geodesic flow on quotients of the hyperbolic plane by a discrete
subgroup of its isometry group. This has been done by Adler and Flatto
\cite{adlerflatto1, adlerflatto3} and by Series \cite{series} for the
continued fraction transformation, and by Adler and Flatto
\cite{adlerflatto2} for R\'enyi's backward continued fraction map (see
also \cite{adler,series2}). In both of these cases, the underlying
surface was the modular surface, which is the quotient of the
hyperbolic plane by the modular group
$\Gamma=PSL(2,\mathbb{Z})$. Representing a system as a factor of a
cross-section of a geodesic flow enables one to derive mechanically
(without guessing) an expression for the invariant density, and to
deduce various properties of the system. Series (\cite{series2},
Problem 5.25(i)) posed the general problem of replicating this idea
for other number- and group-theoretical dynamical systems.

It turns out that the map $T$ also admits such a
representation. Aaronson and Denker \cite{aarodenker} studied a
certain cross-section of the geodesic flow on a different quotient of
the hyperbolic plane, which is the quotient by the congruence subgroup
$\Gamma(2)$ of the modular group, a free group with two
generators. They obtained the map $\hat{T}$ as a factor of that cross
section, and used this to derive results on the asymptotic behavior of
the Poincar\'e series of the group
$\Gamma(\mathbb{C}\setminus\mathbb{Z})$ of deck transformations of
$\mathbb{C}\setminus\mathbb{Z}$. Since in their paper the motivation
was completely different from ours, and the connection to the map $T$
was not known, we find it worthwhile to include here a version of
their result.

{\it \paragraph{Theorem 5.} Let $(\mathbb{H},
(\varphi_t)_{t\in\mathbb{R}})$ be the upper half-plane model of the
hyperbolic plane, with the associated geodesic flow
$\varphi_t:T_1(\mathbb{H})\to T_1(\mathbb{H})$, where
$T_1(\mathbb{H})=\mathbb{H}\times S^1$ is the unit tangent bundle of
$\mathbb{H}$. Define the group of isometries of $\mathbb{H}$
$$ \Gamma(2) = \bigg\{ z\to\frac{az+b}{cz+d} :
\left(\begin{array}{ll}a & b\\c & d\end{array}\right)\in
SL(2,\mathbb{Z}),\ 
\left(\begin{array}{ll}a & b\\c & d\end{array}\right) \equiv
\left(\begin{array}{ll}1 & 0\\0 & 1
\end{array} \right) \textrm{ mod }2 \bigg\}, $$
and let $M = \mathbb{H}/\Gamma(2)$ be the quotient space of
$\mathbb{H}$ by $\Gamma(2)$, which has a fundamental domain
$$ F=\bigg\{ z \in \mathbb{H} : |\textrm{Re }z|<1,\ \left|z\pm
\frac{1}{2}\right| > \frac{1}{2} \bigg\}. $$
Let $\pi:\mathbb{H}\to M$ be the quotient map. Let
$(\overline{\varphi}_t)_{t\in\mathbb{R}}$ be the geodesic flow on
$M$. Let $X' \subset T_1(\mathbb{H})$,
$$ X' = \bigg\{ (z,u)\in \partial F\times S^1 : z+\epsilon u \in
F\textrm{ for small }\epsilon>0 \bigg\}, $$
and let $X\subset T_1(M)$, $X=d\pi(X')$ be the natural section of $M$
corresponding to the fundamental domain $F$ of all inward-pointing
vectors on the boundary of $F$. Let $\tau:X\to X$ be the section- or
first-return map of the geodesic flow, namely
$$ \tau(\omega) = \overline{\varphi}_{t_{\omega}}(\omega), $$ where
$$ t_{\omega} = \inf\{t> 0 : \overline{\varphi}_t(\omega)\in X\}. $$
Then the section map $(X,\tau)$ admits the map $(\qpp,T)$ as a
factor. That is, there exists an (explicit) function $E:X\to\qpp$
such that $T\circ E = E \circ \tau$.}

\bigskip
Theorem 5 is an immediate consequence of Aaronson and Denker's result
(\cite{aarodenker}, section 4) and Theorem 4. To describe explicitly
the factor map $E$, define $p_1(x,\delta,\epsilon)=x$. Then, in the
notation of their paper,
$$E = D \circ p_1 \circ \eta^{-1} \circ \pi_+, $$
with the ``juicy'' parts being our map $D$ and and the map $\pi_+$,
which assigns to a tangent vector the hitting point on the real axis
of the geodesic emanating from the lifting of the tangent vector to
$\partial F\times S^1$. For more details consult \cite{aarodenker}.

The congruence subgroup $\Gamma(2)$ also appears in the paper by
Alperin \cite{alperin}, which discusses the ternary tree of PPTs.

In the next section, we reprove Theorem 1, and show how the linear mappings
involved in the construction of the ternary tree of PPTs can be scaled down
to produce the transformation $T$. This will result in a proof of Theorem
2. In section 3, we prove Theorem 4 and discuss the connection to modified
Euclidean algorithms. The ergodic properties of the system will be derived,
using standard techniques of infinite ergodic theory, proving Theorem 3. In
section 4 we discuss applications to the statistics of expansions of random
points on $\qpp$ and random PPTs. Section 5 has some remarks on points with
special expansions.

\section{Construction of the dynamical system}

\subsection{The piecewise linear transformation}

First, we recall the ideas involved in the proof of Theorem 1. We
follow the elegant exposition of \cite{lonnemo}.

We shall consider solutions of \eqref{eq:pythagoras} with
$\textrm{gcd}(a,b)=1$, and $c>0$. Define
\begin{eqnarray*}
\ppt &=&
  \{ (a,b,c)\in\mathbb{Z}^3 : \textrm{gcd}(a,b)=1,\ \ a,b,c>0,
\ a^2+b^2=c^2 \},
\\ & & \qquad  \textrm{ the set of PPTs, and} \\
\sppt &=&
  \{ (a,b,c)\in\mathbb{Z}^3 : \textrm{gcd}(a,b)=1,\ \ c>0,
\ a^2+b^2=c^2 \},
\\ & & \qquad  \textrm{ the set of \emph{signed} PPTs}.
\end{eqnarray*}
The basic observation is that the equation \eqref{eq:pythagoras} has
three symmetries. Two of them are the obvious symmetries $a\to-a$,
$b\to-b$ (we ignore the symmetry $c\to-c$, since we are only
considering solutions with $c>0$). The third symmetry is not so
obvious, but becomes obvious when the correct change of variables is
applied. Define new variables $m,n,q$ by
$$ \begin{array}{lll} m &=& c-a \\ n &=& c-b \\ q &=& a+b-c
   \end{array}
 \quad \longleftrightarrow \quad
\begin{array}{lll} a &=& q+m \\ b&=&q+n \\ c&=&q+m+n \end{array}
$$
In the new variables, \eqref{eq:pythagoras} becomes
\begin{equation}\label{eqref:pythagorasmnq}
q^2 = 2 m n.
\end{equation}
There is therefore a third natural involution on the set $\sppt$ of
solutions of \eqref{eq:pythagoras}, given in $m,n,q$ coordinates by
$q\to-q$. So, starting from a solution $(a,b,c)\in\sppt$ with
associated variables $(m,n,q)$ and setting $q'=-q, m'=m,n'=n$, we
arrive at a new solution $(a',b',c')$ given by
$$ \begin{array}{lllllll}
a' &=& q'+m' &=& a-2q &=& 2c-a-2b \\
b' &=& q'+n' &=& b-2q &=& 2c-2a-b \\
c' &=& q'+m'+n' &=& c-2q &=& 3c-2a-2b,
\end{array} $$
or in matrix notation
$$ \left( \begin{array}{l}a'\\b'\\c'\end{array}\right) =
\left(\begin{array}{lll}-1&-2&2\\-2&-1&2\\-2&-2&3\end{array}\right)
\left( \begin{array}{l}a\\b\\c\end{array}\right) =: I
\left( \begin{array}{l}a\\b\\c\end{array}\right).
$$
The matrix $I$ is an involution, i.e. $I^2=id_3$. It is easy to see
that $a',b'$ are also coprime, and
$$ c' = 3c-2(a+b) \ge 3c - 2\sqrt{2}\cdot \sqrt{a^2+b^2} =
(3-2\sqrt{2})c > 0. $$
So $I$ maps $\sppt$ to itself. 
%Now, $\sppt$ may be decomposed as
%$$ \sppt = \ppt \cup \sppt_{(+-)} \cup \sppt_{(-+)} \cup \sppt_{(--)}
%\cup \{ (\pm 1,0,1), (0,\pm 1,1) \}, $$
%where $\sppt_{(\epsilon\delta)} = \{ (a,b,c)\in \sppt :
%\sgn(a)=\epsilon, \sgn(b)=\delta \},\ \ \epsilon,\delta\in\{+,-\}$.
$I$ has the fixed points $(1,0,1)$ and $(0,1,1)$. We claim that
$I(\ppt) = \sppt\setminus(\ppt\cup \{(1,0,1),(0,1,1)\})$. Indeed, this
simply means that if $(a,b,c)\in\ppt$, then at least one of $a',b'$ is
negative, or in other words that $2c < \max(a+2b,2a+b)$. Assume for
concreteness that $a>b$, then $2c < 2a+b$ if $2<2x+y$, where
$(x,y)=(a/c,b/c) \in \qpp\cap\{x>y\}$, or equivalently if
\begin{equation}\label{eq:scalarprod}
\frac{2}{\sqrt{5}} < \big\langle\ (x,y),\ (2/\sqrt{5},1/\sqrt{5})
\ \big\rangle. \end{equation}
The point $(2/\sqrt{5},1/\sqrt{5})$ lies on the arc $\qpp\cap\{x>y\}$,
and one checks easily that there is an equality at one end $(1,0)$ of
the arc, and a strict inequality at the other end
$(1/\sqrt{2},1/\sqrt{2})$. So \eqref{eq:scalarprod} holds.

Having shown that if $(a,b,c)\in\ppt$, then $(a',b',c')$ is a signed
PPT with one of $a,b$ negative, we can forget about the signs of
$a',b'$ and obtain a new triple $(a'',b'',c'')=(|a'|,|b'|,c')$. $c'$
is strictly less than $c$, since $c'=c-2 q$ and
$q=a+b-c=a+b-\sqrt{a^2+b^2}>0$ on $\ppt$. The new triple will be a
PPT, except when $(a,b,c)=(3,4,5)$ or $(4,3,5)$, in which case
$(a'',b'',c'')$ will equal $(1,0,1)$ or $(0,1,1)$, respectively. If
$(a'',b'',c'')$ is a PPT, there are precisely three PPTs $(a,b,c)$
leading to it via this procedure -- corresponding to the three possible
sign patterns $a'<0<b';\ \ a',b'<0;\ \ a'>0>b'$ (we ruled out
$a',b'>0$) -- and they can easily be recovered, as follows: If
$a'<0<b'$, then
\begin{multline*}
\qquad
\left( \begin{array}{l}a''\\b''\\c''\end{array}\right) =
\left( \begin{array}{r}-a'\\b'\\c'\end{array}\right)=
\left( \begin{array}{ccc}-1&0&0\\0&1&0\\0&0&1\end{array}\right)
I \left( \begin{array}{l}a\\b\\c\end{array}\right) \\
\implies
\left( \begin{array}{l}a\\b\\c\end{array}\right) =
I \left( \begin{array}{ccc}-1&0&0\\0&1&0\\0&0&1\end{array}\right)
\left( \begin{array}{l}a''\\b''\\c''\end{array}\right) =
M_1 \left( \begin{array}{l}a''\\b''\\c''\end{array}\right). \qquad
\qquad
\end{multline*}
(with $M_1$ as in \eqref{eq:matrices}). Similarly we get
$$
\left( \begin{array}{l}a\\b\\c\end{array}\right) =
M_2 \left( \begin{array}{l}a''\\b''\\c''\end{array}\right),
\textrm{ if $a',b'<0$, or }
\left( \begin{array}{l}a\\b\\c\end{array}\right) =
M_3 \left( \begin{array}{l}a''\\b''\\c''\end{array}\right),
\textrm{ if $a'>0>b'$}.
$$
We are now ready to prove Theorem 1. First, from the above discussion
it follows that $M_1,M_2,M_3$ take PPTs to PPTs with a strictly larger
third coordinate. In particular, any triple of one of the forms
\eqref{eq:expansion1}, \eqref{eq:expansion2} is a PPT. Next, for
$(a,b,c)\in\ppt$ define
\begin{eqnarray*}
S(a,b,c) &=& (|2-a-2b|,|2-2a-b|,3-2a-2b), \\
\delta(a,b,c) &=& \left\{ 
      \begin{array}{ll} 1 & 2c-a-2b<0< 2c-2a-b \\
                 2 & 2c-a-2b, 2c-2a-b<0 \\
                 3\qquad & 2c-a-2b>0>2c-2a-b \end{array}\right.,\\
\delta_k(a,b,c) &=& \delta(S^{k-1}(a,b,c)), \quad
k=1,2,\ldots,n(a,b,c), \\
n(a,b,c) &=& \max\{ n\ge 0 : S^n(a,b,c) \in \ppt \}.
\end{eqnarray*}
The above discussion can be summarized by
the equations
\begin{equation}\label{eq:deltam}
\left( \begin{array}{l}a\\b\\c \end{array}\right) =
M_{\delta(a,b,c)} (S(a,b,c))^{\textrm{t}},
\qquad
\delta\left( M_d
\left( \begin{array}{l}a\\b\\c \end{array}\right) \right) = d\ \ \ 
(d=1,2,3).
\end{equation}
Let $(a,b,c)\in\ppt$ with $a$ odd and $b$ even. It is easy to see that
$M_1,M_2,M_3$ preserve the parity of $a,b$, so $(a,b,c)$ cannot have a
representation \eqref{eq:expansion2}. We claim that it satisfies
\eqref{eq:expansion1} with $d_k = \delta_k(a,b,c)$,\ \
$k=1,\ldots,n(a,b,c)$, and that this representation is unique. The
proof is by induction on $c$. The claim holds for the basic triple
$(3,4,5)$, because $M_1,M_2,M_3$ increase the third coordinate. Assume
that it holds for all odd-even PPTs with third coordinate $<c$. Then
in particular this is true for $(a',b',c')=S(a,b,c)$, since we know
that $c'<c$. So we may write
$$
\left( \begin{array}{l}a'\\b'\\c' \end{array}\right) =
M_{e_1}M_{e_2}\ldots M_{e_m}
\left(\begin{array}{l}3\\4\\5\end{array}\right)
$$
with $ e_k = \delta_k(a',b',c')$,\ \ $1\le k\le
m=n(a',b',c')=n(a,b,c)-1$.  We have
$$ e_k = \delta(S^{k-1}(a',b',c')) = \delta(S^k(a,b,c)) =
\delta_{k+1}(a,b,c)= d_{k+1}, $$
where we denote $d_k = \delta_k(a,b,c)$. Therefore, by \eqref{eq:deltam},
$$ \left( \begin{array}{l}a\\b\\c \end{array}\right) =
M_{\delta(a,b,c)}
\left( \begin{array}{l}a'\\b'\\c' \end{array}\right) =
M_{d_1} M_{d_2} \ldots M_{d_{m+1}}
\left( \begin{array}{l}3\\4\\5 \end{array}\right),
$$
which is our claimed representation. Uniqueness follows immediately by
noting that by \eqref{eq:deltam}, $d_1$ is determined by the sign
pattern of $(2c-a-2b,2c-2a-b)$, and continuing by induction. Theorem 1
is proved.

\subsection{Scaling the transformation}

It is now easy to rescale the transformation $S$ to obtain a
transformation $T$ from $\qpp$ to its closure.
If $(x,y)=(a/c,b/c)\in\qpp\cap{\mathbb Q}^2$ is a rational
point of $\qpp$, which corresponds to the PPT $(a,b,c)$, then the
triple $S(a,b,c)=(|2c-a-2b|,|2c-2a-b|,3c-2a-2b)$ corresponds to the
point 
$$ \left( \frac{|2c-a-2b|}{3c-2a-2b},\frac{|2c-2a-b|}{3c-2a-2b}
\right) = \left( \frac{|2-x-2y|}{3-2x-2y},\frac{|2c-2a-b|}{3c-2a-2b}
\right) $$
in $\overline{\qpp}$. This precisely accounts for our definition of
$T$ in Theorem 2. To explain why the function $d$ is the correct
rescaling of $\delta$, observe, for example, that $\delta(a,b,c) = 1$
iff
$$ 2-x-2y < 0 < 2-2x-y \iff
\begin{array}{l}
\frac{2}{\sqrt{5}} < 
\big\langle\  (x,y), (1/\sqrt{5},2/\sqrt{5})\ \big\rangle \\ \ \\
\big\langle\  (x,y), (2/\sqrt{5},1/\sqrt{5})\ \big\rangle <
\frac{2}{\sqrt{5}}, \end{array}
$$
which some inspection reveals to hold exactly on the (open) circular
arc connecting the point $(1,0)$ with the point $(4/5,3/5)$. This
corresponds to the condition $x/y>4/3$ in the definition of
$d$. Similarly, it can be checked that $\delta(a,b,c)=2$ if $(x,y)$
lies on the circular arc between $(4/5,3/5)$ and $(3/5,4/5)$, and
$\delta(a,b,c)=3$ if $(x,y)$ lies on the circular arc between
$(3/5,4/5)$ and $(0,1)$. These arcs form the \emph{generating
partition} of the ternary expansion -- see Figure 2.

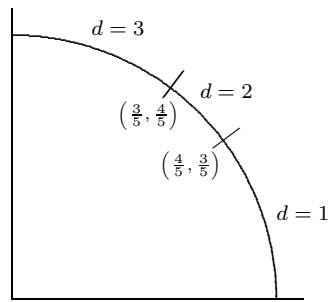
\begin{figure}[h!]
\begin{center}
\begin{picture}(110,110)(0,0)
\arc(100,0){90}
\put(0,0){\line(1,0){110}}
\put(0,0){\line(0,1){110}}
\put(57,76){\line(3,4){8}}
\put(76,57){\line(4,3){10}}
\put(40,68){\tiny$\left(\frac35,\frac45\right)$}
\put(56,49){\tiny$\left(\frac45,\frac35\right)$}
\put(100,30){\scriptsize$d=1$}
\put(71,76){\scriptsize$d=2$}
\put(30,100){\scriptsize$d=3$}
\end{picture}
\caption{The quadrant $\qpp$ and the generating partition}
\end{center}
\end{figure}

Figure 3 shows the graph of the map obtained by parametrizing the
quadrant $\qpp$ in terms of the angle (multiplied by $2/\pi$, to
obtain a map on the interval $(0,1)$). As Theorem 4 may imply, this is not
the best parametrization, but it gives a good graphical illustration
of the behavior of the mapping $T$. Note that, contrary to appearance,
the map is not linear on the middle interval!

\begin{figure}
\begin{center}
\includegraphics{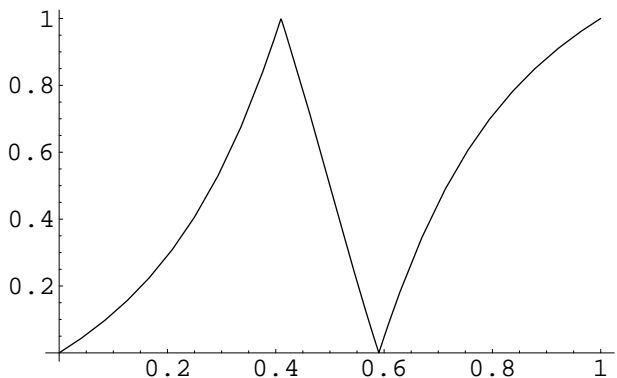}
\caption{The conjugate map $F^{-1}\circ T\circ F$, where
$F(t)=(\cos(\pi t/2),\sin(\pi t/2))$}
\end{center}
\end{figure}

We have proved part (i) of Theorem 2. For part (ii), observe that all
the iterates of an irrational point $(x,y)\in\qpp$ under $T$ remain
irrational, since $T$ is defined by piecewise rational functions with
integer coefficients which are invertible (so as before, given
$d(x,y)$, we can recover $(x,y)$ from $T(x,y)$ by a rational function
with integer coefficients, hence if $T(x,y)$ is rational, so is
$(x,y)$). For $T^n(x,y)$ to be in $\overline{\qpp}\setminus\qpp =
\{(1,0),(0,1)\} $, we must have $T^{n-1}(x,y)=(3/5,4/5)$ or
$(4/5,3/5)$, and this cannot happen for an irrational point. Therefore
$T^n(x,y)$ is defined for all $n\ge 1$, as claimed. The resulting
sequence of digits $d_k = d(T^{k-1}(x,y))$ cannot terminate with an
infinite succession of 1's. Indeed, on the first interval
$(0,\frac{2}{\pi}\arctan(3/4))$ of the generating partition of the
mapping $G=F^{-1}\circ T\circ F$ (Figure 3) it is easy to verify that
$G'$ is increasing and satisfies $G'(0)=1$. Therefore if $t$ is some
point on this interval, then the sequence of iterates $G^k(t)$
satisfies $G^k(t) \ge G'(t)^k\cdot t$ and therefore must eventually
leave this interval (so the corresponding succession of 1's in the
expansion will terminate). A symmetrical argument applies for the
third interval of the generating partition, implying that no infinite
expansion terminates with an infinite succession of 3's, and part (ii)
is proved.

Turn to the final part (iii) of Theorem 2. Again we use the mapping
$G$ in Figure 3. Let $I_1=(0,\frac{2}{\pi}\arctan(3/4))$,
$I_2=(\frac{2}{\pi}\arctan(3/4), \frac{2}{\pi}\arctan(4/3))$,
$I_3=(\frac{2}{\pi}\arctan(4/3),1)$ be the intervals of the generating
partition of $G$.  Given an infinite expansion $(d_k)_{k\ge 1}$ that
does not terminate with an infinite succession of $1$'s or an infinite
succession of $3$'s, we must show that there is a unique number
$t\in(0,1)$ such that $G^{k-1}(t) \in I_{d_k}$ for all $k\ge 1$.

Consider, for any $n\ge 1$, the \emph{cylinder set}
$$ A_n = \{ t\in(0,1) : G^{k-1}(t)\in I_{d_k}\textrm{ for
}1\le k\le n \}. $$
$(A_n)_{n\ge 1}$ is a decreasing sequence of non-empty open
intervals. By compactness, the intersection of their closures contains
at least one point $t$. The condition on the sequence $(d_k)_{k\ge 1}$
implies that $t$ is in fact in the intersection of the \emph{open}
intervals; otherwise, $t$ is an endpoint, say of $A_n$, but that would
imply that $d_k = 1$ for all $k>n$ or $d_k=3$ for all $k>n$.

We have shown existence of a number with a prescribed expansion. But
uniqueness also follows, since, as Figure 3 shows, any appearance of a
``2'' digit, or a non-1 digit following a succession of 1's, or a
non-3 digit following a succession of 3's, entails a shrinkage of the
corresponding set $A_n$ by at least a constant factor bounded away
from 1. So the diameter of the $A_n$'s goes to 0, and their
intersection contains at most one point. Theorem 2 is proved.

(Here is another argument demonstrating uniqueness: any two irrational
points on $\qpp$ are separated by a rational point; after a finite
number of applications of $T$, the rational point will be mapped to
$(3/5,4/5)$ or to $(4/5,3/5)$, and the images of the two irrational
points will be contained in different elements of the generating
partition, implying a different first digit in their expansions.)

\section{The modified Euclidean algorithm}

\subsection{Some computations}

The inverse function of $D$ is easily computed to be
$$ D^{-1}(x,y) = \frac{1-x}{y}. $$
Using this, a routine computation, which we omit, shows that indeed
$$\hat{T}=(D^{-1}\circ T\circ D). $$
We show that the measure $\nu$ is $\hat{T}$-invariant. If
$d\nu(t)=f(t)dt$, the invariant density must satisfy
\begin{equation}\label{eq:invariant}
f(t) = \sum_{u=\hat{T}^{-1}(t)}
f(u)\cdot\frac{1}{|\hat{T}'(u)|}.
\end{equation}
The inverse branches of $\hat{T}$ are given by
\begin{multline}\label{eq:inversebranches}
\qquad\qquad\qquad\qquad\qquad\ \ \ 
F_1(t) = \frac{t}{1+2t} \in (0,1/3), \\
\qquad F_2(t) = \frac{1}{2+t} \in (1/3,1/2), \\
F_3(t) = \frac{1}{2-t} \in (1/2,1), \qquad\qquad\qquad\qquad\ \ 
\end{multline}
for which
\begin{eqnarray*}
|T'(F_1(t))| &=& (1-2F_1(t))^{-2} = (1+2t)^2, \\
|T'(F_2(t))| &=& F_2(t)^{-2} = (2+t)^2, \\
|T'(F_3(t))| &=& F_3(t)^{-2} = (2-t)^2.
\end{eqnarray*}
So \eqref{eq:invariant} reduces to
\begin{equation}\label{eq:invariant2}
f(t) = f\left(\frac{t}{1+t}\right)\frac{1}{(1+2t)^2} +
 f\left(\frac{1}{2+t}\right)\frac{1}{(2+t)^2} +
 f\left(\frac{1}{2-t}\right)\frac{1}{(2-t)^2}.
\end{equation}
Check directly that $f(t)=1/(t(1-t))$ satisfies \eqref{eq:invariant2}.

To complete the proof of Theorem 4, we need to verify that $\mu$ is
the push-forward of the measure $\nu$ under $D$. Denote
$$ x = x(t) = \frac{1-t^2}{1+t^2}, \qquad y = y(t) =
\frac{2t}{1+t^2}. $$
Compute:
$$ \frac{1}{\sqrt{(1-x)(1-y)}} = \left( \frac{2t^2}{1+t^2}\cdot
\frac{(1-t)^2}{1+t^2} \right)^{-1/2} = \frac{1+t^2}{\sqrt{2}\,
t(1-t)}. $$
\begin{eqnarray*}
ds &=& \sqrt{dx^2+dy^2} = \sqrt{x'(t)^2+y'(t)^2} dt \\ &=&
\sqrt{ \left(\frac{-4t}{(1+t^2)^2}\right)^2 +
       \left(\frac{2(1-t^2)}{(1+t^2)^2}\right)^2 } dt =
\frac{2\,dt}{1+t^2}.
\end{eqnarray*}
$$ \implies
\frac{dt}{\sqrt{2}\,t(1-t)} = \frac{ds}{\sqrt{(1-x)(1-y)}}, $$
as claimed.

Note that this also proves that $\mu$ is $T$-invariant. This fact
could be checked directly, of course.

\subsection{Interpretation as a Euclidean algorithm}

The ordinary Euclidean algorithm takes a pair of positive integers
$(x,y)$ with $x>y$ and returns the pair $(y,x\textrm{ mod }y)$. After
a finite number of iterations of this mapping, $y$ will be equal to
$0$ and $x$ will be equal to the g.c.d. of the original pair.

Many variants of this algorithm have been analyzed, where various
alternatives to simple division with remainder are used. The study of such
algorithms, related of course to continued fraction variants, is a huge
subject which it is beyond the scope of this paper to describe. See
\cite{schweiger}; sections 4.5.2-4.5.3 in \cite{knuth} and the references
there; and \cite{baladivallee} for some more recent developments.

The standard Euclidean algorithm has a more ancient version, known as
the \emph{slow}, or \emph{subtractive} Euclidean algorithm, where
subtraction is used instead of division, so $(x,y)$ are mapped to
$(\max(x-y,y),\min(x-y,y))$. Clearly the standard algorithm is nothing
more than a speeding-up of this algorithm. One may scale by always
replacing $x$ by 1 and $y$ by the ratio $y/x$. This leads to the
interval map $R:(0,1)\to(0,1)$,

\vspace{-20.0 pt}
$$
R(t) = \left\{ \begin{array}{ll} \frac{t}{1-t} & 0 < t \le 1/2, \\
  \frac{1-t}{t} & 1/2 < t < 1 \end{array}\right. =
\begin{array}{l} \\
\includegraphics{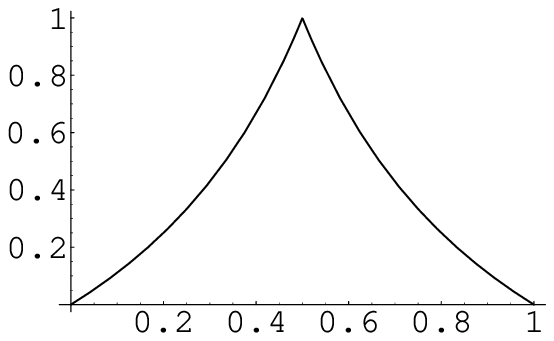}
\end{array}
$$
We now observe that the map $\hat{T}$ is itself the scaling of a
modified algorithm, defined by the mapping
\begin{eqnarray*} (x,y),\ \ x>y &\longrightarrow& \begin{array}{ll}
 (x-2y, y) & \textrm{if }x-2y > y, \\
 (y,x-2y)  & \textrm{if }y \ge x-2y > 0, \\
 (y,2y-x)  & \textrm{if }x-2y \le 0. \end{array}
\end{eqnarray*} 
Here is a sample execution sequence of this algorithm:
$$ (155,100) \longrightarrow (100,45) \longrightarrow (45,10)
\longrightarrow (25,10) \longrightarrow (10,5) \longrightarrow (5,0).
$$
This algorithm can be used to compute g.c.d.'s, just like its famous
kin: The last output which differs from the one preceding it is of either
the form $(a,a)$ or $(a,0)$, where $a$ is the g.c.d. of the two original
integers.

\subsection{Ergodic properties of $\hat{T}$}

To prove Theorem 3, we study the somewhat simpler measure preserving
system $((0,1),\hat{T},\nu)$. Since this is now represented as an
interval map, we can use standard techniques of ergodic theory.

Define $\hat{d}(t) = d(D(t))$. Define the set $J=(1/5,2/3)$. An
alternative description for $J$ is as the cylinder set
\begin{multline*} \quad
J = \bigg\{ t\in(0,1) : \hat{d}(t)=2\textrm{ or }\big(\hat{d}(t)=1
\textrm{ and }\hat{d}(\hat{T}(t))\neq 1\big) \\ \textrm{ or }
\big(\hat{d}(t)=3\textrm{ and }\hat{d}(\hat{T}(t))\neq 3\big) \bigg\}.
\qquad\qquad
\end{multline*}
By Theorem 2(ii), for any irrational $t\in(0,1)$, $\hat{T}^n(t) \in J$
for some $n\ge 0$. In other words,
$$ (0,1) = \bigcup_{n=0}^\infty \hat{T}^{-n}(J)\ \ \textrm{a.e.} $$
This implies that $\hat{T}$ is conservative, by \cite{aaroerg}, Theorem
1.1.7.

To prove that $\hat{T}$ is ergodic, we pass to the \emph{induced system}
$ (J, \hat{T}_J, \mu_{|J}) $, where
\begin{eqnarray*}
\hat{T}_J(t) &=& \hat{T}^{\varphi_J(t)}(t), \\
\varphi_J(t) &=& \inf\{ n\ge 1 : T^n(t)\in J \}.
\end{eqnarray*}
By the explicit description of $\hat{T}_J$ given in \cite{aarodenker}, p.
16, it follows (\cite{aarodenker}, Lemma 5.2) that $\hat{T}_J:J\to J$ is a
topologically mixing Markov map which is uniformly expanding with bounded
distortion, i.e., satisfies
$$
\inf_{t\in J} |\hat{T}_J'(t)| > 1, \qquad
\sup_{t\in J} \frac{|\hat{T}_J''(t)|}{(\hat{T}_J'(t))^2} < \infty.
$$
Therefore (\cite{aaroerg}, Theorem 4.4.7) it is exact, and in particular
it is ergodic. It follows (\cite{aaroerg}, Proposition 1.5.2(2)) that
$\hat{T}$ is itself ergodic. This completes the proof of Theorem 3.

\section{Expansions of random $\qpp$-points and random PPTs}

\subsection{Random points on $\qpp$}

The invariant measure $\mu$ becomes infinite near the two ends of the
quadrant $\qpp$. This means that in a typical expansion, the digits
``1'' and ``3'' will occur infinitely more often than the middle digit
``2''. However, for any two digit sequences, even ones that contain
the digits ``1'' and ``3'', we can ask about their relative density of
occurence in the expansion of typical points.

{\it
\paragraph{Theorem 6.}
Let $I_1=(0,1/3), I_2=(1/3,1/2), I_3=(1/2,1)$ be the intervals of the
generating partition of $\hat{T}$.
For $(d_1,\ldots,d_n)\in\cup_{\ell=1}^\infty\{1,2,3\}^\ell$, define
$$ A(d_1,\ldots,d_n) = \nu\bigg( \bigcap_{j=1}^n 
\hat{T}^{-j+1}(I_{d_j}) \bigg) =
\nu\bigg( \big(F_{d_1}\circ F_{d_2} \circ \ldots \circ F_{d_n}\big)
\big((0,1)\big)
\bigg), $$
with $F_1, F_2, F_3$ as in \eqref{eq:inversebranches}.
Let $(d_1,\ldots,d_n), (e_1,\ldots,e_m) \in \cup_{\ell=1}^\infty
\{1,2,3\}^\ell$. For almost every $(x,y)\in\qpp$, the limit
$$ \lim_{N\to\infty} \frac{\# \bigg\{ 0 \le k \le N :
 d(T^{k+j-1}(x,y)) = d_j,\ \ 1\le j\le n \bigg\} }
{\# \bigg\{ 0 \le k \le N :
 d(T^{k+j-1}(x,y)) = e_j,\ \ 1\le j\le m \bigg\} } $$
exists and is equal to $A(d_1,\ldots,d_n)/A(e_1,\ldots,e_m)$.
}

\paragraph{Proof.} This is an immediate consequence of Hopf's ergodic
theorem, applied to the two indicator functions of the cylinder sets
$\cap_{j=1}^n \hat{T}^{-j+1}(I_{d_j})$ and
$\cap_{j=1}^m \hat{T}^{-j+1}(I_{e_j})$.
\qed

\paragraph{Example.} An easy computation gives
$$ \frac{A(1,2)}{A(1,3)} = \frac{\log(4/3)}{\log(3/2)} \approx
\frac{0.2876}{0.4055}. $$
Therefore, in a typical expansion, when a run of consecutive 1's
breaks, the next digit will be a ``2'' with probability
$0.2876/(0.2876+0.4055) \approx 0.415$, or a ``3'' with probability
$0.4055/(0.2876+0.4055) \approx 0.585$.

\subsection{Random PPTs}

PPTs have finite expansions and form a subset of $\qpp$ of measure 0. So,
as the analogous studies of continued fraction expansions of rational
numbers (a.k.a. analysis of the Euclidean algorithm) have shown, analyzing
their behavior can be significantly more difficult than the behavior of
expansions of random points on $\qpp$. We outline here a technique for
easily deducing some of the properties of the expansion by relating the
discrete model to the continuous one. We mention some open problems which
may be approachable using more sophisticated methods such as those used in
\cite{baladivallee}, and which we hope to address at a later date.

Our model for random PPTs will be the discrete probability space
$$ \ppt_N = \{ (a,b,c)\in \ppt : c \le N \}, $$
equipped with the uniform probability measure
$\mathbb{P}_N$. Analogous results can easily be formulated, using the
same ideas presented here, for other natural models, e.g., a uniform
choice of $(a,b,c)\in PPT$ with $|a|,|b|\le N$.

We discuss the distribution of the individual digits in the expansion. Let
$\lambda$ be the uniform arc-length measure on $\qpp$, normalized as a
probability measure. That is, $d\lambda = (2/\pi)ds$. We need the following
simple lemma.

{\it
\paragraph{Lemma 7.} Under the measure $\mathbb{P}_N$, the random vector
$(a/c,b/c)$ converges in distribution to $\lambda$.
}

\paragraph{Proof.} Observe the following fact from elementary number
theory: the coprime pairs $(m,n)$ such that $m,n$ are of opposite
parity have a local density of $4/\pi^2$ in the lattice
$\mathbb{Z}^2$, in the following sense: for any bounded open set
$D\subset \mathbb{R}^2$, we have
\begin{multline}\label{eq:density}
\frac{1}{x^2}\#\bigg\{ (m,n)\in\mathbb{Z}^2 : 
\left(\frac{m}{x},\frac{n}{x}\right)\in D,\
\textrm{gcd}(m,n)=1,\ \ m+n\equiv 1 (\textrm{mod }2) \bigg\} \\
\xrightarrow[x\to\infty]{}
\frac{4}{\pi^2}\textrm{area}(D). \qquad\qquad\qquad\qquad\qquad
\end{multline}
First, this is true when $D$ is a rectangle $(0,A)\times(0,B)$. To
prove this, define for $i=0,1$,
$$ g_i(u; k) = \#\bigg\{ j\in\mathbb{Z} : 0 < j < u,\ \ j\equiv i
(\textrm{mod 2}),\ \ k\ |\ j \bigg\}. $$
Let $(\mu_2(k))_{k\ge 1}$ be the coefficients of the Dirichlet series
$$ \beta(s):=
\prod_{p>2\textrm{ prime}} (1-p^{-s})^{-1} = \sum_{k=1}^\infty
\mu_2(k) k^{-s}. $$
Then the left-hand side of \eqref{eq:density} is equal to
$$ \frac{1}{x^2}\sum_{k=1}^\infty \mu_2(k) \big[
g_0(Ax; k)g_1(Bx; k)+g_1(Ax; k)g_0(Bx; k) \big]. $$
Since clearly $|g_i(u; k)-(u/2k)|\le 1$ (for $k$ odd), this is easily
seen to be $c A B + O((\log x)/x)$ as $x\to\infty$, where
$$ c = \frac{1}{2} \sum_{k=1}^\infty \mu_2(k)k^{-2} =
\frac{1}{2}\beta(2) = \frac{1}{2}\cdot \frac{4}{3}\zeta(2) =
\frac{4}{\pi^2}, $$
proving our claim.

It follows, by taking unions and differences, that \eqref{eq:density}
is true for $D$ any finite union of rectangles with sides parallel to
the coordinate axes, and therefore by approximation for any bounded
and open $D$.

For $0< t\le \pi/2$, denote
\begin{eqnarray*}
\textrm{arc}(t) &=& \bigg\{ (\cos u,\sin u) : 0 < u < t \bigg\}, \\
\textrm{sector}(t) &=& \bigg\{ (x,y) : x,y>0,\ \ x^2+y^2 \le 1,\ \
\arctan(y/x) < t \bigg\}.
\end{eqnarray*}
By the parametric solution \eqref{eq:parametric} we have as
$N\to\infty$
\begin{eqnarray*}
\mathbb{P}_N\bigg( (a,b,c)\in \ppt_N : (a/c,b/c)\in\textrm{arc}(t)
\bigg)\qquad\qquad\qquad\qquad
\end{eqnarray*}

\vspace{-15.0 pt}
{\scriptsize
\begin{eqnarray*}
&=& \frac{\# \bigg\{ (m,n)\in\mathbb{Z}^2 :\ \ 
m>n>0,\ \ \gcd(m,n)=1,\ \ m+n\equiv 1 (2),
\ \ m^2+n^2\le N,\ \ 
\arctan\left(\frac{2mn}{m^2-n^2}\right)<t
\bigg\}}{\# \bigg\{ (m,n)\in\mathbb{Z}^2 :\ \ 
m>n>0,\ \ \gcd(m,n)=1,\ \ m+n\equiv 1 (2),
\ \ m^2+n^2\le N \bigg\}} \\
&=& \frac{\# \bigg\{ (m,n)\in\mathbb{Z}^2 :\ \ \textrm{gcd}(m,n)=1,\ \ 
m+n\equiv 1 (2),\ \ 
\left(\frac{m}{\sqrt{N}},\frac{n}{\sqrt{N}}\right) \in
\textrm{sector}(t/2) \bigg\}}
{\# \bigg\{ (m,n)\in\mathbb{Z}^2 :\ \ \textrm{gcd}(m,n)=1,\ \ 
m+n\equiv 1 (2),\ \ 
\left(\frac{m}{\sqrt{N}},\frac{n}{\sqrt{N}}\right) \in
\textrm{sector}(\pi/4) \bigg\}} 
\end{eqnarray*}}

\vspace{-30.0pt}
\begin{eqnarray*}
&=& (1+o(1)) \frac{4\pi^{-2}\cdot\textrm{area}(\textrm{sector}(t/2))N}
{4\pi^{-2}\cdot\textrm{area}(\textrm{sector}(\pi/4))N} =
(1+o(1))\frac{2t}{\pi}.
\qquad\qquad\qquad\qquad
\end{eqnarray*}
This is exactly the claim of the Lemma.
\qed

\bigskip
Define the \emph{Perron-Frobenius operator} of $T$ as the operator
$H:L_1(\qpp,\lambda)\to L_1(\qpp,\lambda)$,
\begin{eqnarray*}
(Hf)(x,y) &=&
\frac{1}{3+2x-2y}\cdot
f\left(\frac{2+x-2y}{3+2x-2y},\frac{2+2x-y}{3+2x-2y}\right) \\ & & +
\frac{1}{3+2x+2y}\cdot
f\left(\frac{2+x+2y}{3+2x+2y},\frac{2+2x+y}{3+2x+2y}\right) \\ & & +
\frac{1}{3-2x+2y}\cdot
f\left(\frac{2-x+2y}{3-2x+2y},\frac{2-2x+y}{3-2x+2y}\right).
\end{eqnarray*}
$H$ is also known as the \emph{transfer operator} corresponding to
$T$. It has the property that if the random vector $(X,Y)$ on $\qpp$
has distribution $f(x,y)d\lambda$, then $T(X,Y)$ has distribution
$(Hf)(x,y)d\lambda$. We skip the simple computation that verifies this
claim.

{\it \paragraph{Theorem 8.} Under $\mathbb{P}_N$, the distribution of
$d(T^{n-1}(a/c,b/c))$, the $n$-th digit in the expansion of a random
PPT $(a,b,c)\in\ppt_N$, converges to the distribution of $d(x,y)$
under the measure
$$ d\lambda_n(x,y) = (H^{n-1}({\bf 1}))(x,y) d\lambda(x,y), $$
where ${\bf 1}$ is the constant function $1$.
}

\paragraph{Proof.} This is immediate from Lemma 7 and the definition
of $H$.
\qed

\bigskip
Theorem 8 answers the question of the limiting distribution as $N\to\infty$
of the digits $d_n$ in the expansion of a random PPT; however, it does not
give good asymptotic information on the behavior of these distributions as
$n\to\infty$. In fact, this is not a very interesting question: since the
invariant measure is infinite, the density $H^n({\bf 1})$ will become for
large $n$ more and more concentrated around the singular points $(1,0)$ and
$(0,1)$ (it is possible to make this statement more precise, but we do not
pursue this slightly technical issue here).

Here's one way to amend the situation in a way that enables formulating
interesting quantitative statements concerning the distribution of the
digits, which we mention briefly without going into detail: replace the
expansion $(d_k)_{k=1}^n$ by a new expansion $(e_j)_{j=1}^\ell$, by
dividing the $(d_k)$ into blocks consisting of 1's and 3's and terminating
with a 2; so for instance, the expansion $(1,1,2,2,3,1,3,2,3,3,3,1,2,1)$
will be replaced by $(112,2,3132,33312,1)$.  The new expansion corresponds
to the induced transformation $T_B$, where $B$ is the middle arc in the
generating partition. Cylinder sets of $T_B$ can be easily computed. The
invariant measure is the restriction $\mu_{|B}$, a \emph{finite} measure, so
normalize it to be a probability measure. $T_B$ is easily shown to be exact
as in section 3.3. A theorem analogous to Theorem 8 above can be proved, to
the effect that the $n$-th digit in the ``new'' expansion of a random PPT
converges in distribution to the distribution of $d_{\textrm{new}}(x,y)$
(the first ``new'' digit) under the measure whose density with respect to
$\lambda$ is the $(n-1)$-th iterate of the Perron-Frobenius operator of
$T_B$ applied to the constant function ${\bf 1}$.  Since $T_B$ is mixing,
these densities will actually converge to the invariant density. So after
each occurence of a ``2'' in the original expansion, there are well-defined
statistics for the sequence of digits that follows up to the next ``2''.

We conclude with some open problems: study the expectation, the variance
and the limiting distribution of the \emph{length} of the expansion of a
random element of $PPT_N$, as $N\to\infty$. Generalize to arbitrary
cost-functions, as in \cite{baladivallee}.

\section{Some special expansions}

From a number-theoretic standpoint, it is interesting to study points on
$\qpp$ with special expansions. As the examples in section 1 show, simple
periodic expansions seem to correspond to simple quadratic points on
$\qpp$. It is easy to see that any eventually-periodic expansion
corresponds to the image under $D$ of a quadratic irrational. Do all
quadratic irrationals have eventually periodic expansions?

We also found empirically the expansions
$$ \begin{array}{llcll}
e^i&=&(\cos 1,\sin 1) &=& [3,1^2,3,1^4,3,1^6,3,1^8,3,\ldots], \\
e^{i/2}&=& (\cos 1/2,\sin 1/2) &=& [1,3,1^5,3,1^9,3,1^{13},3,\ldots].
\end{array}$$
where $1^k$ means a succession of $k$ 1's. The first equation
can be proved by observing
that $D^{-1}(\cos 1,\sin 1)=(1-\cos(1))/\sin(1)=\tan(1/2)$, and that the
approximations
$$ \big(
F_3\circ F_1^2 \circ F_3 \circ F_1^4 \circ \ldots \circ F_3 \circ
F_1^{2k} \big)(1)
$$
($F_1, F_2, F_3$ as in \eqref{eq:inversebranches}) have the continued
fraction expansions
$$
\frac{1}{1}
\begin{array}{ll}\\+\end{array}
\frac{1}{1}
\begin{array}{ll}\\+\end{array}
\frac{1}{4}
\begin{array}{ll}\\+\end{array}
\frac{1}{1}
\begin{array}{ll}\\+\end{array}
\frac{1}{8}
\begin{array}{ll}\\+\end{array}
\frac{1}{1}
\begin{array}{ll}\\+\end{array}
\frac{1}{12}
\begin{array}{ll}\\+\end{array}
\ldots
\begin{array}{ll}\\+\end{array}
\frac{1}{4k}. $$
Thus our expansion reduces to the known (\cite{sloane},
sequence A019425) infinite continued fraction expansion
$$ \tan(1/2) =
\frac{1}{1}
\begin{array}{ll}\\+\end{array}
\frac{1}{1}
\begin{array}{ll}\\+\end{array}
\frac{1}{4}
\begin{array}{ll}\\+\end{array}
\frac{1}{1}
\begin{array}{ll}\\+\end{array}
\frac{1}{8}
\begin{array}{ll}\\+\end{array}
\frac{1}{1}
\begin{array}{ll}\\+\end{array}
\frac{1}{12}
\begin{array}{ll}\\+\end{array}
\ldots
$$
The second equation is proved similarly. It is interesting to wonder
whether other ``nice'' expansions exist for simple points on $\qpp$.

\section*{Acknowledgements}

Thanks to Jon Aaronson and to Boaz Klartag for helpful discussions.


\begin{thebibliography}{02}

\bibitem{aaroerg} J. Aaronson,
\newblock An Introduction to Infinite Ergodic Theory.
\newblock Mathematical Surveys and Monographs 50,
\newblock Amer. Math. Soc, Providence, RI, 1997.

\bibitem{aarodenker} J. Aaronson, M. Denker,
\newblock The Poincar\'e series of $\mathbb{C}\setminus\mathbb{Z}$.
\newblock {\it Ergodic Theory Dynam. Systems} 19 (1999), 1--20.

\bibitem{adler} R. L. Adler,
\newblock Geodesic flows, interval maps, and symbolic dynamics.
\newblock In: Ergodic Theory, Symbolic Dynamics and Hyperbolic
Spaces, 93--123.
\newblock Oxford University Press, Oxford, 1991.

\bibitem{adlerflatto1} R. L. Adler, L. Flatto,
\newblock Cross section maps for geodesic flows I (the modular
surface).
\newblock In: Ergodic theory and dynamical systems,
vol. 2 (College Park, Maryland 1979-1980), 103--161.
\newblock Progr. Math., 21, Birkh\"auser, Boston, Mass., 1982.

\bibitem{adlerflatto2} R. L. Adler, L. Flatto,
\newblock The backward continued fraction map and geodesic flow.
\newblock {\it Ergodic Theory Dynam. Systems} 4 (1984), 487--492.

\bibitem{adlerflatto3} R. L. Adler, L. Flatto,
\newblock Cross section map for the geodesic flow on the modular
surface.
\newblock In: Conference in modern analysis and probability
(New Haven, Conn., 1982), 9--24.
\newblock Contemp. Math., 26, Amer. Math. Soc., Providence, RI, 1984.

\bibitem{alperin} R. Alperin,
\newblock The Modular tree of Pythagoras.
\newblock Preprint,
2000. \texttt{http://www.arxiv.org/abs/math.HO/0010281}.

\bibitem{baladivallee} V. Baladi, B. Valle\'e,
\newblock Euclidean algorithms are Gaussian. Preprint, 2003.
\texttt{http://arxiv.org/abs/cs.DS/0307062}

\bibitem{barning} F. J. M. Barning,
\newblock On Pythagorean and quasi-Pythagorean triangles and a
generation process with the help of unimodular matrices. (Dutch)
\newblock {\it Math. Centrum Amsterdam Afd. Zuivere Wisk.}, ZW-011 (1963).

\bibitem{billingsley} P. Billingsley,
\newblock Ergodic Theory and Information.
\newblock John Wiley \& Sons, New York-London-Sydney 1965.

\bibitem{gollnick} J. Gollnick,  H. Scheid, J. Z\"ollner,
\newblock Rekursive Erzeugung der primitiven pythagoreischen
Tripel. (German)
\newblock {\it Math. Semesterber.} 39 (1992), 85--88.

\bibitem{hall} A. Hall,
\newblock Genealogy of Pythagorean triads.
\newblock {\it Math. Gazette} 54:390 (1970), 377--379.

\bibitem{hardywright} G. H. Hardy, E. M. Wright,
\newblock An Introduction to the Theory of Numbers, 5th ed.
\newblock Oxford University Press, Oxford, 1985.

\bibitem{jaeger} J. Jaeger,
\newblock Pythagorean number sets. (Danish.)
\newblock {\it Nordisk Mat. Tidskr.} 24 (1976), 56--60, 75.

\bibitem{kanga} A. R. Kanga,
\newblock The family tree of Pythagorean triples.
\newblock {\it Bull. Inst. Math. Appl.} 26 (1990), 15--17.

\bibitem{knuth}
D. E. Knuth,
\newblock The Art of Computer Programming, vol. 2: Seminumerical
Algorithms, 3rd. ed.
\newblock Addison-Wesley, 1998.

\bibitem{kristensen} E. Kristensen,
\newblock Pythagorean number sets and orthonormal matrices. (Danish)
\newblock {\it Nordisk Mat. Tidskr.} 24 (1976), 111--122, 135.

\bibitem{lonnemo} A. L\"onnemo,
\newblock The trinary tree underlying primitive pythagorean triples.
\newblock In: Cut the Knot, Interactive Mathematics Miscellany and
Puzzles. Alex Bogomolny (Ed.), \\
\texttt{http://www.cut-the-knot.org/pythagoras/PT\_matrix.shtml}.

\bibitem{mccullough} D. McCullough,
\newblock Height and excess of Pythagorean triples.
\newblock Preprint.

\bibitem{preau} P. Pr\'eau,
\newblock Un graphe ternaire associ\'e \`a l'\'equation
$X\sp 2+Y\sp 2=Z\sp 2$.
\newblock {\it C. R. Acad. Sci. Paris Ser. I Math.} 319 (1994),
665--668.

\bibitem{schweiger} F. Schweiger,
\newblock Ergodic Theory of Fibred Systems and Metric Number Theory.
\newblock Clarendon Press, Oxford, 1995.

\bibitem{series} C. Series,
\newblock The modular surface and continued fractions.
\newblock {\it J. London Math. Soc.} (2) 31 (1985), 69--80.

\bibitem{series2} C. Series,
\newblock Geometrical methods of symbolic coding.
\newblock In: Ergodic Theory, Symbolic Dynamics and Hyperbolic
Spaces, 125--151.
\newblock Oxford University Press, Oxford, 1991.

\bibitem{sloane} N. J. A. Sloane, editor (2003),
\newblock The On-Line Encyclopedia of Integer Sequences.
\newblock \texttt{http://www.research.att.com/~njas/sequences/}.

\end{thebibliography}
\end{document}